\documentclass[12pt]{article}

\newtheorem{definition}{{\bf Definition}}[section]
\newtheorem{theorem}[definition]{{\bf Theorem}}
\newtheorem{lemma}[definition]{{\bf Lemma}}
\newtheorem{proposition}[definition]{{\bf Proposition}}

\newtheorem{example}[definition]{{\bf Example}}

\def\ffi{\varphi}
\def\<{\langle}
\def\>{\rangle}

\def\iH{{\cal H}}

\def\iA{{\cal A}}
\def\iB{{\cal B}}
\def\iC{{\cal C}}

\def\Tr{{\rm Tr}}

\usepackage{amscd,amsmath,amssymb,graphicx}

\begin{document}

\begin{center}
 {\Large Quasi-orthogonal subalgebras of matrix algebras}@
\end{center}

\begin{center}
Hiromichi Ohno \\ 
\end{center}

\section{Introduction}

In the theory of quantum mechanics, 
an $n$-level system is described by the algebra $M_n = M_n({\mathbb C})$
of $n \times n$ complex matrices.
The matrix algebra of a composite system consisting of an $n$-level and an
$m$-level system is $M_n \otimes M_m \simeq M_{nm}$.
A subalgebra of $M_k$ corresponds to a subsystem of a $k$-level quantum system.
In this paper, subalgebras contain the identity and are closed under the adjoint
operation of matrices, that is, they are unital $*$-subalgebras.
The algebra $M_n$ can be endowed by the inner product
$\langle A,B \rangle = {\rm Tr}(A^*B)$ and it becomes a Hilbert space.
Two subalgebras ${\cal A}_1$ and ${\cal A}_2$ are called quasi-orthogonal if
${\cal A}_1 \ominus {\mathbb C}I \perp {\cal A}_2 \ominus {\mathbb C}I$.

We consider pairwise quasi-orthogonal subalgebras
$\iA_1, \iA_2, \ldots , \iA_l$ in $M_{p^{kn}}$ which are isomorphic to $M_{p^k}$ for $k\ge 1$, $n\ge 2$ and  a prime number
$p$  with $p \ge 3$.
The aim of this paper is to obtain the maximum $l$.
The case $p=2$, $n=2$ and $k=1$ is shown in \cite{OhPeSt, PDcomp, nofive} and
the maximum is $4$.

The motivations of this problem are followings.
If a total system $M_n \otimes M_m$ has a statistical operator $\rho$,
we can reconstruct the reduced density $\rho_n^{(1)} = \Tr_{m} \rho$ in
the subsystem $M_n$,
where $\Tr_{m}$ is a partial trace onto $M_n$.
In order to get more information, we change the density 
$\rho$ by an interaction.
For a Hamiltonian $H$, the new state is
\[
e^{iH} \rho e^{-iH} = W_1 \rho W_1^*
\]
after the interaction. The new reduced density is 
$\rho_n^{(2)} = \Tr_m W_1 \rho W_1^*$.
By using other interactions, we have a sequence of reduced states 
$\rho_n^{(1)}, \rho_n^{(2)}, \ldots , \rho_n^{(k)}$.
We want to determine the minimum $k$ such that this sequence of 
reduced densities determines $\rho$.

Another reason comes from the relation to the mutually unbiased bases (MUB)
 problem. 
Given a orthonormal basis of an $n$-dimensional Hilbert space $\iH$,
the linear operators diagonal in this basis form a maximal Abelian subalgebra
of $M_n \simeq \iB(\iH)$.
Conversely if $|e_i\>\<e_i|$ are minimal projections 
in a maximal Abelian subalgebra,
then $(|e_i\>)_i$ is a orthonormal basis.
Two maximal Abelian subalgebras of $M_n$ are quasi-orthogonal if and only if
the two corresponding bases $(\xi_i)$ and $(\zeta_j)$ are mutually unbiased
(see Proposition \ref{prop2.1}), that is,
\[
|\< \xi_i , \zeta_j \> |^2 = {1\over n} \qquad (1\le i,j \le n).
\]
Mutually unbiased bases are very interesting from many point of view and
the maximal number of such bases is not known for arbitrary $n$.
We want to study mutually unbiased bases in terms of 
quasi-orthogonal subalgebras (see \cite{BBTV, BSTW, PiRu}).

In Section \ref{sect2}, we consider quasi-orthogonal subalgebras
in $M_{p^n}$ which are isomorphic to $M_p$. In Section \ref{sect3},
quasi-orthogonal subalgebras in $M_{p^{kn}}$ which are 
isomorphic to $M_{p^k}$ are investigated.


\section{Quasi-orthogonal subalgebras in $M_{p^n}$}\label{sect2}

${\cal A}$ is a finite dimensional $C^*$-algebra with the usual trace ${\rm Tr}$ and 
is considered as a Hilbert space under the inner product
\[
\langle A , B \rangle = {\rm Tr}(A^*B)
\]
for any $A,B \in {\cal A}$.

\begin{definition}{\rm
Two subalgebras ${\cal A}_1$ and ${\cal A}_2$ of ${\cal A}$ are called {\it quasi-orthogonal} if 
\[
{\cal A}_1 \ominus {\mathbb C}I \perp {\cal A}_2 \ominus {\mathbb C}I.
\]}
\end{definition}

The equivalent conditions of this definition are followings \cite{Popa}:

\noindent
(a) For any $A_1 \in {\cal A}_1$ and $A_2 \in {\cal A}_2$,
\[
{\rm Tr}(A_1A_2) = { {\rm Tr}(A_1){\rm Tr}(A_2) \over {\rm Tr}(I)}.
\]

\noindent
(b) For any $A_1 \in {\cal A}_1$ and $A_2 \in {\cal A}_2$ with ${\rm Tr}(A_1) = {\rm Tr}(A_2)=0$,
\[
{\rm Tr}(A_1A_2) =0.
\]

The theory of quasi-orthogonal subalgebras is related to the theory of mutually unbiased bases.
This is stated as follows:

\begin{proposition}\label{prop2.1}
The maximal number of mutually unbiased bases in ${\mathbb C}^d$ is equal to
the maximal number of pairwise quasi-orthogonal subalgebras in $M_d$ which are
isomorphic to ${\mathbb C}^d$.
\end{proposition}
\begin{proof}
For orthonormal bases $\{\xi_i\}_{1\le i \le d}$ and 
$\{\zeta_j\}_{1\le j \le d}$ of ${\mathbb C}^d$,
define minimal projections in $M_d$ by
$P_i = |\xi_i \rangle \langle \xi_i |$ and
$Q_j =|\zeta_j \rangle \langle \zeta_j |$.
Since ${\rm Tr} (P_i Q_j) = |\langle \xi_i  , \zeta_j \rangle|^2$, the condition
\begin{eqnarray*}
|\langle \xi_i  , \zeta_j \rangle|^2 = {1\over d} \qquad (1\le i,j \le d)
\end{eqnarray*}
is equivalent to the condition
\[
{\rm Tr} (P_i Q_j) = {1\over d} = {\Tr (P_i) \Tr (Q_j) \over \Tr (I)}
\qquad (1\le i,j \le d).
\]
From (a), we conclude that $\{\xi_i\}$ and $\{ \zeta_j \}$ are
 mutually unbiased bases
if and only if the algebras generated by $\{P_i\}$ and $\{Q_j\}$ are quasi-orthogonal.
This prove the assertion.
\end{proof}

It is known that there exist $d+1$ mutually unbiased bases in ${\mathbb C}^d$
if $d$ is a power of a prime number.
Therefore there are $d+1$ quasi-orthogonal subalgebras in $M_d$ which are isomorphic to
${\mathbb C}^d$ if $d$ is a power of a prime number.

Now we consider the quasi-orthogonal subalgebras in 
$M_{p^2} =M_p \otimes M_p$ which are isomorphic to $M_p$,
where $p$ is a prime number with $p \ge 3$.
We will show that we can construct $p^2 +1$ pairwise quasi-orthogonal subalgebras.

Define the unitary operators $W$ and $S$ in $M_p$ by
\[
W=
\left[\begin{array}{ccccc}
1 & 0 &0& \cdots & 0 \\
0 & \lambda &0 & \cdots & 0 \\
0 & 0& \lambda^2  & \cdots & 0 \\
\vdots & \vdots &\vdots & \ddots & \vdots \\
0 & 0 &0&  \cdots & \lambda^{p-1}
\end{array}\right], \qquad
S=
\left[\begin{array}{ccccc}
0 & 0 &  \cdots & 0 & 1 \\
1 & 0 &  \cdots & 0 &0 \\
0 & 1 &  \cdots & 0 &0\\
\vdots & \vdots & \ddots &\vdots & \vdots \\
0 & 0  & \cdots &1 & 0
\end{array}\right],
\]
where $\lambda = e^{2\pi i / p}$.
Then these two unitary operators have following properties:

\begin{enumerate}
\item[i.] $S^p = W^p = I$.
\item[ii.] The set $\{ S^i W^j\}_{0\le i,j \le p-1}$ is a natural orthogonal basis of $M_p$.
\item[iii.] Since $SW = \lambda^{-1} WS$, we have
\begin{equation}\label{eq4.1}
S^{k_1} W^{l_1} S^{k_2} W^{l_2} = \lambda^{k_2l_1 } S^{k_1 + k_2} W^{l_1 + l_2}.
\end{equation}
\item[iv.]  From (\ref{eq4.1}), $S^{k_1} W^{l_1}$ and $S^{k_2} W^{l_2}$ commute if and only if
$k_1l_2 =k_2 l_1$ mod $p$.
\end{enumerate}

We consider the commutativity condition iv 
in the context of a vector space over
 $Z_p$, where $Z_p$ is a finite field with $p$ elements.
Let $Z_p^4 = \{ (k_1, l_1, k_2, l_2) \, | \, k_1,l_1,k_2,l_2 \in Z_p\}$ 
be a vector space over $Z_p$ and define
 a natural homomorphism $\pi$ (up to scalar multiple)
from $Z_p^4 $ to $M_p \otimes M_p$ by
\[
\pi (k_1,l_1,k_2,l_2) = S^{k_1}W^{l_1} \otimes S^{k_2} W^{l_2}.
\]
We denote a {\it symplectic product} by
\[
u \circ u' = k_1l_1' - k_1' l_1  + k_2 l_2' - k_2' l_2  \qquad {\rm mod} \,\,p,
\]
where $u = (k_1,l_1,k_2,l_2)$ and $u'= (k_1',l_1',k_2',l_2')$. 
From (\ref{eq4.1}), 
\begin{equation}\label{eq4.2}
\pi(u) \pi(u') = \lambda^{-u\circ u'} \pi(u')\pi(u).
\end{equation}
Hence $\pi(u) $ and $\pi(u')$ commute
if and only if their symplectic product equals zero.

\begin{lemma}\label{lemma4.1}
If $\pi(u)$ and $\pi(u')$ are not commutative 
for $u = (k_1,l_1,k_2,l_2), \, u'= (k_1',l_1',k_2',l_2') \in Z_p^4$, then
the algebra ${\cal A}$ generated by $\pi(u)$ and $\pi(u')$ is isomorphic to $M_p$.
\end{lemma}
\begin{proof}
From the assumption, $u\circ u' \neq 0$.
Define a map $\rho$ from $\{S, W^{u\circ u'} \}$ to ${\cal A}$ by
\[
\rho(S) = \pi(u), \qquad \rho(W^{u\circ u'}) = \pi(u').
\]
From (\ref{eq4.2}) and $SW^{u\circ u'} = \lambda^{-u\circ u'} W^{u\circ u'} S$,
the commutativity condition of $\pi(u)$, $\pi(u')$ and that of $S$, $W^{u\circ u'}$ are same.
Therefore $\rho$ can be extended to an isomorphism from 
$M_p$ generated by $S$ and $W^{u\circ u'}$ to ${\cal A}$.
\end{proof}

From this lemma, we need to find such $u$ and $u'$.
Let $D$ be a non-zero integer in $Z^p$ with the requirement
that $D \neq k^2$ mod $p$ for all $k$ in $Z_p$, i.e., $D$ is not a quadratic residue of $p$.
For any $a_0, a_1 \in Z_p$, we define a subspace of $Z_p^4$ by
\begin{eqnarray*}
C_{a_0, a_1} = \{ b_0 (1, a_1,0,a_0) +b_1 (0,a_0, -1 ,a_1D) \,|\, b_0,b_1 \in Z_p \},
\end{eqnarray*}
where scalar multiplication and addition are defined by a natural way.
Moreover  put
\begin{eqnarray*}
C_{\infty} = \{ b_0 (0, 1,0,0) +b_1 (0,0, 0 ,1) \,|\, b_0,b_1 \in Z_p \}.
\end{eqnarray*}


\begin{lemma}\label{lemma2.3}
The only vector common to any pair of above subspaces is $(0,0,0,0)$.
In particular, the subspaces expect $(0,0,0,0)$ partition $Z_p^4 \backslash \{ (0,0,0,0)\}$.
\end{lemma}
\begin{proof}
Since there are $p^2+1$ subspaces and each subspace has $p^2$ elements,
it is enough to prove that the intersection of any two subspaces is $\{(0,0,0,0)\}$.
It is easy to see $C_{a_0,a_1} \cap C_\infty = \{(0,0,0,0)\}$. 
Therefore we prove that $C_{a_0,a_1} \cap C_{a_0',a_1'} =\{(0,0,0,0)\}$ if $a_0 \neq a_0'$ or
$a_1 \neq a_1'$. 

Assume $b_0 (1, a_1,0,a_0) +b_1 (0,a_0, -1 ,a_1D) =b_0' (1, a_1',0,a_0') +b_1' (0,a_0', -1 ,a_1'D) $,
then from the first and third components we have $b_0 =b_0'$ and $b_1 = b_1'$.
Similarly, from the second and fourth components we are led to the equations
\begin{eqnarray*}
a_1 b_0 + a_0 b_1 &=& a_1' b_0 + a_0' b_1 \\
a_0 b_0 + a_1 b_1 D &=& a_0'b_0 + a_1'b_1D.
\end{eqnarray*}
These equations can be rewritten as a matrix equation
\[
\left[ 
\begin{array}{cc}
b_1 & b_0 \\
b_0 & b_1 D
\end{array}
\right]
\left[ 
\begin{array}{c}
a_0 -a_0' \\
a_1-a_1'
\end{array}
\right]
=
\left[ 
\begin{array}{c}
0 \\
0
\end{array}
\right].
\]
If $b_0 = b_1 =0$, then the common element is $(0,0,0,0)$.
Therefore we assume $b_0\neq 0$ or $b_1 \neq 0$.
Then the above matrix is invertible, indeed
\[
\left[ 
\begin{array}{cc}
b_1 & b_0 \\
b_0 & b_1 D
\end{array}
\right]^{-1}
=
(b_1^2 D - b_0^2)^{-1}
\left[ 
\begin{array}{cc}
b_1D & -b_0 \\
-b_0 & b_1 
\end{array}
\right].
\]
Here we use that $b_1^2 D \neq b_0^2$ mod $p$ from the assumption of $D$.
This implies $a_0 = a_0'$ and $a_1 = a_1'$ which is a contradiction.
\end{proof}

Since $(1, a_1,0,a_0)\circ (0,a_0, -1 ,a_1D) =2a_0$, if $a_0 \neq 0$ then
the algebra ${\rm span} \{\pi (C_{a_0,a_1})\} $ generated by $\pi(1,a_1,0,a_0)$
and $\pi(0,a_0,-1,a_1D)$ is isomorphic to $M_p$
by Lemma \ref{lemma4.1}. But if $a_0 =0$, then the algebra is commutative and hence
${\rm span} \{ \pi(C_{a_0,a_1}) \} \simeq {\mathbb C}^{p^2}$.
Therefore we need to mix the subspaces with  $a_0 =0$.

For any $a \in Z_p$, define subspaces by
\begin{eqnarray*}
D_a &=& \{ b_0(1,1,-a,aD) + b_1 (1,2, -a, 2aD) \,|\,b_0 ,b_1 \in Z_p \},\\
 D_\infty &=& \{ b_0 (0,0,1,0) + b_1(0,0,0,1) \,|\, b_0,b_1\in Z_p \}.
\end{eqnarray*}

\begin{lemma}\label{lemma2.4}
The only vector common to any pair of above subspaces is $(0,0,0,0)$.
Moreover we have
\[
\bigcup_{a\in Z_p} D_a \cup D_\infty = \bigcup_{a_1 \in Z_p} C_{0,a_1} \cup C_\infty.
\]
\end{lemma}
\begin{proof}
The proof of the first assertion is same as Lemma \ref{lemma2.3}.
 To show the second assertion, it is enough to prove
$D_a, D_\infty \subset \bigcup_{a_1 \in Z_p} C_{0,a_1} \cup C_\infty$.
Indeed, since each subspaces has $p^2$ elements and the 
intersection is trivial,
the numbers of elements in both sets are equal.

First consider a element $(0,0,b_0,b_1)$ in $D_\infty$.
If $b_0 = 0$, then $(0,0,0,b_1) \in C_\infty$. If $b_0 \neq 0$, then
\[
(0,0,b_0,b_1) = -b_0 (0,0,-1, -b_0^{-1}b_1D^{-1} D) \in C_{0, -b_0^{-1}b_1D^{-1}}.
\]
Hence $D_\infty \subset \bigcup_{a_1 \in Z_p} C_{0,a_1} \cup C_\infty$.
Next consider the element $b_0(1,1,-a,aD) + b_1 (1,2, -a, 2aD)$ in $D_a$.
If $b_0 + b_1 = 0$, then 
\[
b_0(1,1,-a,aD) + b_1 (1,2, -a, 2aD) = (0,b_0 +2 b_1,0,ab_0D + 2ab_1D) \in
C_\infty.
\]
If $b_0+b_1 \neq 0$, then 
\begin{eqnarray*}
&& b_0(1,1,-a,aD) + b_1 (1,2,- a, 2aD) \\
&=& (b_0+b_1)\left(1, (b_0+b_1)^{-1}(b_0+2b_1),0,0\right)\\ && + a(b_0+b_1)\left(0,0,-1,(b_0+b_1)^{-1}(b_0+2b_1)D\right)\\
&\in& C_{0, (b_0+b_1)^{-1}(b_0+2b_1)}.
\end{eqnarray*}
Therefore we obtain $D_a \subset \bigcup_{a_1 \in Z_p} C_{0,a_1} \cup C_\infty$.
\end{proof}

Since $(1,1,-a,aD) \circ  (1,2, -a, 2aD) = 1 -a^2D \neq 0$ by the assumption of $D$
and $(0,0,1,0)\circ (0,0,0,1) =1$, we obtain
\begin{eqnarray*}
{\rm span} \{ \pi(D_a) \} \simeq M_p , \\
{\rm span} \{ \pi(D_\infty) \} \simeq M_p
\end{eqnarray*}
by Lemma \ref{lemma4.1}. Consequently we have the next theorem.

\begin{theorem}\label{theorem2.4}
There are $p^2 +1$ pairwise quasi-orthogonal subalgebras in
 $M_{p^2}$ which are isomorphic to
$M_p$.
\end{theorem}

\begin{example}{\rm
Consider a $9$-level quantum system $M_3 \otimes M_3$.
We list $10$ pairwise quasi-orthogonal subalgebras in $M_9$ which are
isomorphic to $M_3$:

\begin{eqnarray*}
 {\rm span}\{\pi(C_{1,0})\}
&=& {\rm span} \{ I\otimes I, S \otimes W, S^2 \otimes W^2,
W \otimes S^2, SW \otimes S^2 W, \\ && S^2 W \otimes S^2W^2,
W^2 \otimes S , SW^2 \otimes SW, S^2W^2 \otimes SW^2 \},\\
{\rm span}\{\pi(C_{1,1})\}
&=& {\rm span} \{ I\otimes I, SW \otimes W, S^2 W^2 \otimes W^2,
W \otimes S^2W^2 , SW^2 \otimes S^2 , \\ && S^2  \otimes S^2W,
W^2 \otimes SW , S \otimes SW^2, S^2W \otimes S\},\\
{\rm span}\{\pi(C_{1,2})\}
&=& {\rm span} \{ I\otimes I, SW^2 \otimes W, S^2W \otimes W^2,
W \otimes S^2W , S \otimes S^2 W^2, \\ && S^2 W^2 \otimes S^2,
W^2 \otimes SW^2 , SW \otimes S, S^2  \otimes SW \}, \\
{\rm span}\{\pi(C_{2,0})\}
&=& {\rm span} \{ I\otimes I, S \otimes W^2, S^2 \otimes W,
W^2 \otimes S^2 , SW^2 \otimes S^2W^2 , \\ && S^2W^2  \otimes S^2W,
W \otimes S , SW \otimes SW^2, S^2W \otimes SW \},\\
{\rm span}\{\pi(C_{2,1})\}
&=& {\rm span} \{ I\otimes I, SW \otimes W^2, S^2 W^2 \otimes W,
W^2 \otimes S^2W^2 , S \otimes S^2W , \\ && S^2W  \otimes S^2,
W \otimes SW , SW^2 \otimes S, S^2 \otimes SW^2 \},\\
{\rm span}\{\pi(C_{2,2})\}
&=& {\rm span} \{ I\otimes I, SW^2 \otimes W^2, S^2 W \otimes W,
W^2 \otimes S^2W, SW \otimes S^2 , \\ && S^2  \otimes S^2W^2,
W \otimes SW^2 , S \otimes SW, S^2W^2 \otimes S\},\\
{\rm span}\{\pi(D_0)\}
&=& {\rm span} \{ I\otimes I, SW \otimes I, S^2 W^2 \otimes I,
SW^2 \otimes I , S^2 \otimes I , \\ && W  \otimes I,
S^2W \otimes I , W^2 \otimes I, S \otimes I\},\\
{\rm span}\{\pi(D_{1})\}
&=& {\rm span} \{ I\otimes I, SW \otimes S^2W^2, S^2 W^2 \otimes SW,
SW^2 \otimes S^2W , S^2 \otimes S , \\ && W  \otimes W^2,
S^2W \otimes SW^2 , W^2 \otimes W, S \otimes S^2\},\\
{\rm span}\{\pi(D_{2})\}
&=& {\rm span} \{ I\otimes I, SW \otimes SW, S^2 W^2 \otimes S^2W^2,
SW^2 \otimes SW^2 , S^2 \otimes S^2 , \\ && W  \otimes W,
S^2W \otimes S^2W , W^2 \otimes W^2, S \otimes S\},\\
{\rm span}\{\pi(D_\infty)\}
&=& {\rm span} \{ I\otimes I, I \otimes S, I \otimes S^2,
I \otimes W , I \otimes SW , \\ && I  \otimes S^2W,
I \otimes W^2 , I \otimes SW^2, I \otimes S^2W^2\},
\end{eqnarray*}
where we denote $D=2$.
}
\end{example}


Next, we consider pairwise quasi-orthogonal subalgebras in $M_{p^n}$ which
are isomorphic to $M_p$.
The dimensions of $M_{p^n}\ominus {\mathbb C} I$ and $M_p \ominus {\mathbb C}I$ are
$p^{2n}-1$ and $p^2-1$, respectively.
Hence such subalgebras are at most $N_n := p^{2n}-1 / p^2-1$.
We construct $N_n$ pairwise quasi-orthogonal subalgebras.

\begin{theorem}
 $M_{p^{n}}$ contains $N_n$ pairwise quasi-orthogonal subalgebras which
are isomorphic to $M_p$.
\end{theorem}
\begin{proof}
The case $n=1$ is trivial and the case $n=2$ is already proven in Theorem \ref{theorem2.4}.
Assume $M_{p^{n-2}}$ contains $N_{n-2}$ pairwise 
quasi-orthogonal subalgebras which
are isomorphic to $M_p$.
Let $\iA \simeq M_p \subset M_{p^{n-2}}$ be one of them.
Then we can assume that $\iA$ is generated by $W$ and $S$.
On the other hand, $M_{p^2}$ contains $p^2 +1$ quasi-orthogonal subalgebras which are isomorphic
to ${\mathbb C}^{2p} \simeq {\mathbb C}^p \otimes {\mathbb C}^p$ from Proposition \ref{prop2.1}
and the theory of mutually unbiased bases.
Let $\iC \simeq {\mathbb C}^p \otimes {\mathbb C}^p \subset M_{p^2}$ be one of them.
Suppose $\iC$ is generated by $W\otimes I$ and $I\otimes W$.

We define a homomorphism $\hat\pi$ 
(up to scalar multiple) from ${Z}_p^4$ to $\iA \otimes \iC$ by
\[
\hat\pi(k, l, m_1, m_2) = S^kW^l \otimes W^{m_1} \otimes W^{m_2}.
\]
For any $a_0,a_1 \in {\mathbb Z}_p$, let
\[
\hat{C}_{a_0,a_1} = \{ b_0 (1,0,a_0,a_1) + b_1(0,1,a_1D,a_0) \,|\, b_0, b_1 \in {\mathbb Z}_p\}
\]
and
\[
\hat{C}_\infty = \{ b_0 (0,0,1,0) + b_1(0,0,0,1) \,|\, b_0, b_1 \in {\mathbb Z}_p\}.
\]
From a similar method of the proof of Lemma \ref{lemma2.3},
the only vector common to any pair of above subspaces
is $(0,0,0,0)$.
Moreover since $W$ and $S$ satisfy (\ref{eq4.1}) and 
$\iC$ is commutative, we have
\[
\hat\pi(1, 0,a_0,a_1)  \hat\pi(0,1,a_1D,a_0)   
   = \lambda^{-1} \hat\pi(0,1,a_1D,a_0)\hat\pi(1,0,a_0,a_1).
\]
Therefore the algebra ${\rm span}\{\hat\pi(\hat{C}_{a_0,a_1})\}$ generated by 
$\hat\pi(1, 0,a_0,a_1)$ and $ \hat \pi(0,1,a_1D,a_0)$ is isomorphic
to the algebra generated by $S$ and $W$, that is, $M_p$.
Furthermore the algebra ${\rm span}\{\hat\pi(\hat{C}_{\infty})\}$ generated by
$\hat\pi (0,0,1,0)$ and $\hat\pi(0,0,0,1)$ is $I \otimes \iC$.
Consequently $\iA \otimes \iC$ can be decomposed by $\iA \otimes I = {\rm span }\{\hat\pi
(\hat{C}_{0,0})\}$, $I\otimes \iC$ and $p^2-1$ pairwise quasi-orthogonal subalgebras
which are isomorphic to $M_p$. We denote the $p^2-1$ pairwise quasi-orthogonal subalgebras
by $\{\iB_{\iA, \iC}^{k}\}_{k=1}^{p^2-1}$.

Let
$\iA_1, \ldots, \iA_{N_{n-2}}$ be pairwise quasi-orthogonal subalgebras
in $M_{p^{n-2}}$ which are isomorphic to $M_p$ and let $\iC_1, \ldots , \iC_{p^2+1}$ be
pairwise quasi-orthogonal subalgebras in  $M_{p^2}$ which are isomorphic to ${\mathbb C}^{p^2}$. 
For each $\iA_{i_1}, \iA_{i_2} \subset M_{p^{n-2}}$ and $\iC_{j_1}, \iC_{j_2} \subset M_{p^2}$ with
$i_1 \neq i_2$ or $j_1 \neq j_2$,
we have
\[
(\iA_{i_1}\ominus {\mathbb C}I) \otimes (\iC_{j_1} \ominus {\mathbb C}I) \perp 
\iA_{i_2} \otimes \iC_{j_2}.
\]
This implies
 $\iB_{\iA_{i_1}, \iC_{j_1}}^{k_1} \ominus {\mathbb C}I \perp
\iB_{\iA_{i_2}, \iC_{j_2}}^{k_2} \ominus {\mathbb C}I$ if $i_1 \neq i_2$ or $j_1 \neq j_2$ or 
$k_1 \neq k_2$.
Furthermore  $\iB_{\iA_{i_1}, \iC_{j_1}}^{k}$ is quasi-orthogonal with 
$\iA_{i_2} \otimes I$ and $I\otimes \iC_{j_2}$ for all $i_2, j_2$.
In consequence,
$M_{p^n}$ is decomposed by $M_{p^{n-2}} \otimes I_{M_{p^2}} = 
{\rm span}\{\iA_i \otimes I \,|\, i  \} $ and 
$ I_{M_{p^{n-2}}}  \otimes M_{p^2}= {\rm span} \{ I\otimes \iC_j \,|
\, j\} $ and $\{ \iB_{ \iA_i,\iC_j }^{k} \}_{i,j,k}$.
From Theorem \ref{theorem2.4}, 
$M_{p^{2}}$ contains $p^2+1$ pairwise quasi-orthogonal subalgebras which are isomorphic
to $M_p$. Therefore we get
\[
{N_{n-2}} + p^2+1 + (p^2-1)(p^2+1){N_{n-2}} 
={N_n}
\]
pairwise quasi-orthogonal subalgebras in $M_{p^n}$ which are isomorphic to $M_p$.
\end{proof}


\section{Quasi-orthogonal subalgebras in $M_{p^{kn}}$}\label{sect3}

In this section, we consider quasi-orthogonal subalgebras in $M_{p^{kn}}$ which
are isomorphic to $M_{p^k}$.
First we consider the case $n=2$.

$GF(p^k)$ denotes a finite field with $p^k$ elements.
Up to isomorphisms, $GF(p^k)$ is unique and is defined using a polynomial
\[
f(x) = c_0 + \cdots + c_{k-1}x^{k-1} + x^k
\]
that is irreducible over the field ${Z}_p$.
Then we can write
\[
GF(p^k) = \{ a(t) = a_0 + a_1 t + \cdots +a_{k-1} t^{k-1}
\, | \, a_i \in {Z}_p , \, 0 \le i \le k-1 \},
\]
where $t$ satisfies $f(t) = 0$.
Let 
\[
GF(p^k)^4 = \{ a=(a^{(1)}(t),a^{(2)}(t),a^{(3)}(t),a^{(4)}(t))
 \, | \, a^{(i)}(t)\in GF(p^k), \, 1\le i \le 4 \}
\]
and define a symplectic product:
\[
a \circ b = a^{(1)}(t)  b^{(2)}(t) -a^{(2)}(t) b^{(1)}(t)
+a^{(3)}(t) b^{(4)}(t)-a^{(4)}(t)b^{(3)}(t)
\]
for $a=(a^{(1)}(t),a^{(2)}(t),a^{(3)}(t),a^{(4)}(t))$ 
and $b = (b^{(1)}(t),b^{(2)}(t),b^{(3)}(t),b^{(4)}(t))$ in $GF(p^k)^4$.
Similarly,  symplectic product on $Z_{p}^{4k}$ is denoted by
\[
u\circ v = \sum_{i=1}^k u_{i}^{(1)}  v_i^{(2)}  - u_{i}^{(2)} v_i^{(1)}
              +u_{i}^{(3)}  v_i^{(4)}  - u_{i}^{(4)} v_i^{(3)}
\]
for $u = (u_i^{(j)})_{1\le i \le k, 1\le j \le 4}$ and 
$v = (v_i^{(j)})_{1\le i \le k, 1\le j \le 4}$ in $Z_p^{4k} = (Z_p^k)^4$.
From \cite{PiRu},
there exist a linear functional $\varphi$ on $GF(p^k)$ and a linear isomorphism $\pi_1$
from $GF(p^k)^4$ to $ Z_p^{4k}$ such that
\begin{equation}\label{eq5.1}
\ffi(a\circ b) = \pi_1(a) \circ \pi_1(b),
\end{equation}
where we consider $GF(p^k)^4$ and $Z_p^{4k}$ as vector spaces over $Z_p$.
Define a natural homomorphism (up to scalar multiple) 
$\pi_2$ from $Z_p^{4k}$ to 
$M_{p^{2k}}= \bigotimes_{i=1}^{k} M_p \otimes \bigotimes_{i=1}^{k} M_p$ by
\[
\pi_2(u) = \bigotimes_{i=1}^k S^{u_i^{(1)}} W^{u_i^{(2)}} \otimes 
\bigotimes_{i=1}^k S^{u_i^{(3)}} W^{u_i^{(4)}}.
\]
By (\ref{eq4.1}), we have
\[
\pi_2(u)\pi_2(v) = \lambda^{-u\circ v} \pi_2(v) \pi_2(u).
\]

\begin{lemma}\label{lemma10.1}
If $a \circ b \neq 0$ for $a,b \in GF(p^k)^4$, 
then the algebra ${\cal A}$ generated by 
\[
\{ \pi_2 \pi_1(p(t) a),  \pi_2  \pi_1(q(t) b) \, | \, p(t) ,q(t) \in GF(p^k) \}
\]
is isomorphic to $M_{p^k}$,
where $p(t) a$ is defined by 
\[
p(t) (a_1(t),a_2(t),a_3(t),a_4(t)) = 
(p(t) a_1(t), p(t) a_2(t), p(t) a_3(t) , p(t) a_4(t)).
\]
\end{lemma}
\begin{proof}
If $\{ p(t) a \,|\, p(t) \in GF(p^k) \} 
\cap \{q(t) b \,|\, q(t) \in GF(p^k) \} \neq \{ (0,0,0,0) \}$, 
there exists $r(t) \in GF(p^k)$ such that
$a = r(t)b$. This shows $a \circ b = r(t) (b\circ b) =0$.
Hence we can assume that  
$\{ p(t) a \} \cap \{q(t) b \} = \{ (0,0,0,0) \}$.

For $q(t) \in GF(p^k)$, let 
$\psi_{q(t)}$ be a functional on $GF(p^k)$ defined by
\[
\psi_{q(t)}(p(t)) = \pi_1(p(t) a) \circ \pi_1(q(t)b),
\]
where we consider $GF(p^k)$ as a vector space over $Z_p$.
By (\ref{eq5.1}),  we have
\[
\psi_{q(t)}(p(t)) = \ffi( p(t)a \circ q(t)b)=  \ffi ( p(t) q(t) (a \circ b)).
\]
Since $a\circ b\neq 0$,
$\psi_{q^{(1)}(t)} =\psi_{q^{(2)}(t)}$ implies $q^{(1)}(t) = q^{(2)}(t)$.
Therefore we obtain $\{ \psi_{q(t)} \, | \, q(t) \in  GF(p^k)\} 
= GF(p^k)^*$,
where $GF(p^k)^*$ is a dual space of $GF(p^k)$.

Let $\{ p^{(i)} (t) \}_{1\le i \le k }$ be a basis of the vector 
space $GF(p^k)$.
Then there exists a basis $\{ q^{(j)} \}_{1\le j \le k} \subset GF(p^k)$
 such that
\[
\ffi (p^{(i)}(t)a \circ q^{(j)} (t)b) =\psi_{q^{(j)}(t)} (p^{(i)}(t)) = \delta_{ij}.
\]
For $1\le i,j \le k$, let $S_i = I^{\otimes i-1} \otimes S \otimes I^{\otimes k-i}$ and
$W_j = I^{\otimes j-1} \otimes W \otimes I^{\otimes k-j}$ in $M_{p^k} = \bigotimes^k M_p$.
We define a map $\rho$ from 
$\{S_i , W_j \,|\, 1\le i,j \le k\}$ to ${\cal A}$ by
\begin{eqnarray*}
\rho(S_i) &=& \pi_2 \pi_1 ( p^{(i)}(t) a), \\
\rho(W_j)&=& \pi_2 \pi_1 (q^{(j)}(t) b).
\end{eqnarray*}
Then the commutativity condition of $\{ S_i, W_j \}$, that is,
\begin{eqnarray*}
S_i W_j &=& \lambda^{-\delta_{ij}} W_j S_i, \\
S_i S_j  &=& S_j S_i, \\
W_i W_j &=& W_j W_i
\end{eqnarray*}
and the commutativity condition of 
$\{\pi_2 \pi_1 ( p^{(i)}(t) a),\pi_2 \pi_1 (q^{(j)}(t) b)
\,|\, 1\le i,j\le k \}$, that is,
\begin{eqnarray*}
\pi_2  \pi_1 ( p^{(i)}(t) a)  \cdot \pi_2  \pi_1 (q^{(j)}(t) b) &=&
\lambda^{-\pi_1(p^{(i)}(t)a) \circ \pi_1(q^{(j)}(t)b)}  \cdot 
 \pi_2  \pi_1 (q^{(j)}(t) b) \cdot \pi_2  \pi_1 ( p^{(i)}(t) a) \\
&=& \lambda^{-\ffi (p^{(i)}(t) a\circ q^{(j)}(t) b)} \cdot 
 \pi_2  \pi_1 (q^{(j)}(t) b) \cdot \pi_2  \pi_1 ( p^{(i)}(t) a) \\
&=&\lambda^{-\delta_{ij}} \cdot \pi_2  \pi_1 (q^{(j)}(t) b) \cdot \pi_2  \pi_1 ( p^{(i)}(t) a), \\
\pi_2  \pi_1 ( p^{(i)}(t) a) \cdot \pi_2  \pi_1 ( p^{(j)}(t) a) &=&
\pi_2  \pi_1 ( p^{(j)}(t) a) \cdot \pi_2  \pi_1 ( p^{(i)}(t) a), \\
\pi_2  \pi_1 (q^{(i)}(t) b) \cdot \pi_2  \pi_1 (q^{(j)}(t) b) &=&
\pi_2  \pi_1 (q^{(j)}(t) b) \cdot \pi_2  \pi_1 (q^{(i)}(t) b)
\end{eqnarray*}
are same. Hence $\rho$ can be extended to an isomorphism from
$M_{p^k}$ generated by $\{ S_i, W_j\}_{1\le i,j \le k}$ to ${\cal A}$.
\end{proof}

From this lemma, we need to find such $a$ and $b$. Let $D$ be a non-zero element in
$GF(p^k)$ with the requirement that $D\neq p(t)^2$ for all $p(t)$ in $GF(p^k)$.
For any $a(t), b(t) \in GF(p^k)$, define a subspace of $GF(p^k)^4$ by
\[
C_{a(t), b(t)} = \{ p(t) ( 1,b(t), 0, a(t)) + q(t)(0, a(t), -1, b(t)D) \, | \, p(t), q(t) \in GF(p^k)\}.
\]
Moreover, put
\[
C_\infty =\{ p(t) (0,1,0,0) + q(t)(0,0,0,1) \, | \, p(t), q(t) \in GF(p^k)\}.
\]
\begin{lemma}\label{lemma10.2}
The only vector common to any pair of above subspaces is $(0,0,0,0)$.
In particular, the subspaces expect $(0,0,0,0)$ partition $GF(p^k)^4 \backslash \{ (0,0,0,0)\}$.
\end{lemma}
\begin{proof}
This proof is same as Lemma \ref{lemma2.3}.
\end{proof}

Since $(1,b(t),0,a(t)) \circ (0,a(t),-1,b(t)D) = 2a(t)$, if $a(t) \neq 0$ then
\[
{\rm span} \{ \pi_2 \pi_1(C_{a(t),b(t)}) \} \simeq M_{p^k}
\]
by Lemma \ref{lemma10.1}. But if $a(t) =0$, the algebra is commutative and 
hence ${\rm span} \{ \pi_2 \pi_1(C_{a(t),b(t)}) \} 
\simeq {\mathbb C}^{p^{2k}}$.
Therefore we need to mix the subspaces with $a(t)=0$.

For any $a(t) \in GF(p^k)$, define subspaces by
\begin{eqnarray*}
D_{a(t)} &=& \{ p(t) (1,1,-a(t), a(t)D) + q(t)(1,2,-a(t), 2a(t)D) 
\,|\, p(t),q(t) \in GF(p^k) \}  \\
 D_\infty &=& \{ p(t) (0,0,1,0) + q(t) (0,0,0,1) 
\,|\, p(t) , q(t) \in GF(p^k) \}.
\end{eqnarray*}

\begin{lemma}
The only vector common to any pair of above subspaces is $(0,0,0,0)$.
Moreover we have
\[
\bigcup_{a(t)\in GF(p^k)} D_{a(t)} \cup D_\infty
= \bigcup_{b(t) \in GF(p^k)} C_{0,b(t)} \cup C_\infty.
\]
\end{lemma}
\begin{proof}
This proof is same as Lemma \ref{lemma2.4}.
\end{proof}

Since $(1,1,-a(t), a(t)D) \circ (1,2,-a(t), 2a(t)D) = 1-a(t)^2D \neq 0$ by the assumption of 
$D$ and 
$(0,0,1,0) \circ (0,0,0,1) = 1$, we have 
\begin{eqnarray*}
{\rm span} \{ \pi_2  \pi_1 (D_{a(t)}) \}  &\simeq M_{p^k}, \\
{\rm span} \{ \pi_2  \pi_1 (D_\infty) \}  &\simeq M_{p^k}
\end{eqnarray*}
by Lemma \ref{lemma10.1}. Consequently we have the next theorem.

\begin{theorem}\label{thm10.3}
There are $p^{2k} +1$ pairwise quasi-orthogonal subalgebras in $M_{p^{2k}}$ which 
are isomorphic to $M_{p^k}$.
\end{theorem}


Next we consider pairwise quasi-orthogonal subalgebras in $M_{p^{kn}}$ which are isomorphic to
$M_{p^k}$.
The dimensions of $M_{p^{kn}}\ominus {\mathbb C}I$ and 
$M_{p^k} \ominus {\mathbb C}I$
are $p^{2kn}-1$ and $p^{2k}-1$, respectively.
Hence such subalgebras are at most $N_{(k,n)} = p^{2kn}-1 / p^{2k}-1$.
We construct $N_{(k,n)} $ pairwise quasi-orthogonal subalgebras.

\begin{theorem}
$M_{p^{kn}}$ contains $N_{(k,n)} $ pairwise quasi-orthogonal subalgebras which are isomorphic to $M_{p^k}$.
\end{theorem}
\begin{proof}
The case $n=1$ is trivial and the case $n=2$ is already proven in Theorem \ref{thm10.3}.
Assume $M_{p^{k(n-2)}}$ contains $N_{(k,n-2)}$ pairwise quasi-orthogonal subalgebras which are
isomorphic to $M_{p^k}$.
Let ${\cal A} \simeq M_{p^k} \subset M_{p^{k(n-2)}}$ be one of them.
Then we can assume that ${\cal A}$ is generated by 
$\{S_i, W_j\}_{1\le i,j \le k}$.
On the other hand, $M_{p^{2k}}$ contains 
$N_{(k,2)} = p^{2k}+1$ quasi-orthogonal subalgebras
which are isomorphic to ${\mathbb C}^{p^{2k}} \simeq {\mathbb C}^{p^k} \otimes {\mathbb C}^{p^k}$ 
from Proposition \ref{prop2.1} and the theory of mutually unbiased bases.
Let ${\cal C} \simeq {\mathbb C}^{p^k} \otimes {\mathbb C}^{p^k} \subset M_{p^{2k}}$ 
be one of them. Suppose ${\cal C}$ is generated by 
$\{W_i \otimes I, I \otimes W_j \}_{1\le i,j\le k}$.

Define a homomorphism $\hat{\pi}_2$ (up to scalar multiple) from $Z_p^{4k}$ to
${\cal A}\otimes {\cal C}$ by 
\[
\hat{\pi}_2(u_i^{(j)}) = 
\bigotimes_{i=1}^k S^{u_i^{(1)}} W^{u_i^{(2)}} \otimes 
\bigotimes_{i=1}^k W^{u_i^{(3)}} \otimes \bigotimes_{i=1}^k W^{u_i^{(4)}},
\]
for $(u_i^{(j)})_{1\le i\le k, 1\le j \le 4} \in Z_p^{4k}$.
We denote other symplectic products by
\[
a \hat\circ b = a^{(1)}(t)  b^{(2)}(t) -a^{(2)}(t) b^{(1)}(t)
\]
for $a=(a^{(1)}(t),a^{(2)}(t),a^{(3)}(t),a^{(4)}(t))$ 
and $b = (b^{(1)}(t),b^{(2)}(t),b^{(3)}(t),b^{(4)}(t))$ in $GF(p^k)^4$, and
\[
u \hat\circ v = \sum_{i=1}^k u_i^{(1)} v_i^{(2)} - u_i^{(2)} v_i^{(1)}
\]
for $u= (u_i^{(j)})_{1\le i\le k, 1\le j \le 4}$ and $v= (v_i^{(j)})_{1\le i\le k, 1\le j \le 4}$ in
$Z_p^{4k}$.
Since ${\cal C}$ is commutative, we have
\[
\hat{\pi}_2(u)\hat{\pi}_2(v) = \lambda^{-u\hat\circ v} \hat{\pi}_2(v)\hat{\pi}_2(u).
\]
Moreover from \cite{PiRu}, there exists a homomorphism $\hat{\pi}_1$
from  $GF(p^k)^4$ to $Z_p^{4k}$ such that
\[
\ffi(a\hat\circ b) = \hat{\pi}_1 ( a) \hat\circ \hat{\pi}_1 ( b).
\]
Here for elements $a =(a^{(1)}(t),a^{(2)}(t),0,0)$ and 
$b = (0,0,b^{(3)}(t),b^{(4)}(t))$ in $GF(p^k)^4$, we can assume
\begin{eqnarray*}
\hat\pi_1 (a) &=& ((u_i^{(1)}), (u_j^{(2)}),0,0) \\
\hat\pi_1 (b) &=& (0,0, (v_i^{(3)}),(v_j^{(4)}))
\end{eqnarray*}
for some $(u_i^{(1)}),(u_j^{(2)}),(v_i^{(3)}),(v_j^{(4)}) \in Z_p^k$.
Then we can prove that 
the algebra generated by 
$\hat{\pi}_2 \hat{\pi}_1(a)$ and $\hat{\pi}_2  \hat{\pi}_1(b)$ is 
isomorphic to $M_{p^k}$, if $a\hat\circ b \neq 0$, by using a similar proof of
 Lemma \ref{lemma10.1}.

For any $a(t), b(t) \in GF(p^k)$, let
\[
\hat{C}_{a(t),b(t)} = \{ p(t) (1,0,a(t), b(t)) + q(t) (0,1, b(t)D ,a(t)) \,|\, p(t),q(t) \in GF(p^k) \}
\]
and 
\[
\hat{C}_\infty = \{ p(t) (0,0,1,0) + q(t)(0,0,0,1) \,|\,  p(t),q(t) \in GF(p^k) \}.
\]
From the similar method of the proof of Lemma \ref{lemma2.3}, the only vector common
to any pair of above subspaces is $(0,0,0,0)$.
Since $(1,0,a(t),b(t))\circ (0,1,b(t)D,a(t)) =1$, we obtain
\[
{\rm span}\{ \hat{\pi}_2 \hat{\pi}_1 (\hat{C}_{a(t),b(t)})\} \simeq M_{p^k}.
\]
Moreover we have
\[
{\rm span}\{ \hat{\pi}_2 \hat{\pi}_1 (\hat{C}_\infty )\} = I \otimes {\cal C} \simeq {\mathbb C}^{p^{2k}}.
\]
Consequently ${\cal A}\otimes {\cal C}$ can be decomposed by 
${\cal A} \otimes I = {\rm span}\{\hat{\pi}_2 \hat\pi_1 (\hat{C}_{0,0})\}$, 
$I\otimes {\cal C} ={\rm span} 
\{ \hat{\pi}_2 \hat\pi_1 (\hat{C}_{\infty} )\}$ and
$p^{2k}-1$ pairwise quasi-orthogonal subalgebras which are isomorphic to $M_{p^k}$. We denote the $p^{2k}-1$ pairwise quasi-orthogonal subalgebras by
$\{ {\cal B}_{{\cal A}, {\cal C}}^k \}_{k=1}^{p^{2k}-1}$.

Let ${\cal A}_1 , \ldots , {\cal A}_{N_{(k,n-2)}}$ be pairwise quasi-orthogonal subalgebras in
$M_{p^{k(n-2)}}$ which are isomorphic to $M_{p^k}$ and let 
${\cal C}_1, \ldots , {\cal C}_{N_{(k,2)}}$ be pairwise quasi-orthogonal subalgebras in
$M_{p^{2k}}$ which are isomorphic to ${\mathbb C}^{p^{2k}}$.
Since 
\[
({\cal A}_{i_1} \ominus {\mathbb C}I) \otimes ({\cal C}_{j_1} \ominus {\mathbb C}I)
\perp {\cal A}_{i_2} \otimes {\cal C}_{j_2}
\]
for $i_1 \neq i_2$ or $j_1 \neq j_2$,
${\cal B}_{{\cal A}_{i_1},{\cal C}_{j_1}}^{k_1}$ and
${\cal B}_{{\cal A}_{i_2},{\cal C}_{j_2}}^{k_2}$ are quasi-orthogonal
if $i_1 \neq i_2$ or $j_1 \neq j_2$ or $k_1 \neq k_2$.
Furthermore, ${\cal B}_{{\cal A}_{i_1},{\cal C}_{j_1}}^{k_1}$ is
quasi-orthogonal with ${\cal A}_{i_2}\otimes I_{M_{p^{2k}}}$ and 
$I_{M_{p^{k(n-2)}}} \otimes {\cal C}_{j_2}$ 
for any $i_2$ and $j_2$.
Therefore
we can decompose $M_{p^{kn}}$ by
$M_{p^{k(n-2)}} \otimes I_{M_{p^{2k}}} = 
{\rm span} \{ {\cal A}_i \otimes I \}$, 
$I_{M_{p^{k(n-2)}}} \otimes M_{p^{2k}} = 
{\rm span}\{ I \otimes {\cal C}_j \}$ and 
$\{ {\cal B}_{{\cal A}_i, {\cal C}_j}^k \}_{i,j,k}$.
In consequence, we get
\[
N_{(k,n-2)} + N_{(k,2)} + (p^{2k}-1)N_{(k,2)}N_{(k,n-2)} = N_{(k,n)}
\]
pairwise quasi-orthogonal subalgebras in $M_{p^{kn}}$ which are isomorphic to $M_{p^k}$.
\end{proof}


\end{document}